\documentclass[12pt]{article}
\usepackage{amsmath,amssymb}
\usepackage{amsthm}
\usepackage{epsfig}
\newcommand{\Z}{{\mathbb Z}}
\newcommand{\cx}{{\mathbf c}}
\newcommand{\ux}{{\mathbf u}}
\newcommand{\vx}{{\mathbf v}}
\newcommand{\yx}{{\mathbf y}}

\newtheorem{theorem}{Theorem}
\newtheorem{proposition}{Proposition}
\newtheorem{corollary}{Corollary}
\newtheorem{lemma}{Lemma}
\newtheorem{conjecture}{Conjecture}

\begin{document}

\title{\bf $k$-Dependence and Domination in Kings Graphs%
\footnote{or ``Too Many Kings and There Goes the Neighborhood''}}
\author{
\textbf{
 Eugen J. Ionascu$^1$, Dan Pritikin$^2$, Stephen E. Wright$^2$
}\\
$^1$\small Department of Mathematics, Columbus State University\\
\small Columbus, GA 31907; \texttt{ionascu\_eugen@colstate.edu}\\
\small Honorific Member of the Romanian Institute\\
\small of Mathematics 
``Simion Stoilow''\\
$^2$\small Department of Mathematics and Statistics, Miami University\\
\small Oxford, OH 45056; 
  \texttt{ pritikd@muohio.edu, wrightse@muohio.edu}
}

\date{\today\ \textbf{Draft}}

%

\maketitle

\section{Introduction}

Among combinatorial chessboard problems, one of the most basic is 
the following, essentially found in \cite{yaglom}:
\begin{quote}
What is the maximum number of kings that can be placed on an 
$m\times n$ board, so that no two squares occupied by kings share 
a side or a corner (i.e., no king ``attacks'' another)?
\end{quote}
By placing kings in row $i$ and column $j$ when $i$ and $j$ are 
both odd, we see that at least $\lceil m/2 \rceil \lceil n/2 
\rceil$ squares can be occupied by kings.  With some thought, one 
can prove that this is optimal.  For $m$ and $n$ large, one 
learns that the best one can do is place kings on about $1/4$ of 
the squares.  In this paper, we study the following general 
version of this problem:
\begin{quote}
Given a whole number $k$, what is the maximum number $s$ of kings that
can be placed on an $m \times n$ board, so that no king attacks 
more than $k$ other kings? When $m$ and $n$ are large, how large 
can the density $s/(mn)$ be?
\end{quote}
For most choices of $k=0,\ldots,8$, there is a tidy solution: an 
upper bound can be proved by a short elementary argument, and an 
arrangement of kings can be constructed to show that the upper 
bound is tight.  These limiting densities are given in section 
\ref{upperbound}.  However, tight upper bounds are not yet known 
for either $k=4$ or $k=5$.  It is easy to construct arrangements 
of kings (on arbitrarily large boards) that achieve the densities 
of $3/5$ and $9/13$ for $k=4$ and $5$, respectively.  We 
conjecture that these are indeed the maximum limiting densities.

The story in the present article concerns the struggle in 
supporting this conjecture by good upper bounds, as well as the 
variety of rival techniques used for different values of $k$.  
Along the way, we make elementary use of graph theory, number 
theory, group theory, real analysis, and integer linear 
programming.  We believe the methods of the present paper can 
provide the basis for undergraduate research projects on related 
problems.


\section{Notation and Terminology}

We have already deviated from traditional chess in several ways:  
the board's length and width are arbitrary; each chess piece is a 
king with no associated color; we are concerned with optimal 
arrangements of pieces, rather than actual chess moves.  We 
actually go a few steps further.  First, we generalize the 
discussion to address the density problem of placing kings on 
multidimensional chessboards.  Second, it is useful to also treat 
\emph{toroidal} boards allowing ``wrap-around"; these provide an 
idealization with the same limiting densities as non-toroidal 
boards, but with a simpler analysis.  Third, some results are 
stated in terms of arbitrary graphs.  These three extensions also 
serve to identify possible areas for undergraduate research.

We adopt notation and terminology from graph theory by referring 
to board squares as \emph{vertices}.  We let $K[n_1,\ldots,n_d]$ 
denote the $n_1 \times \cdots \times n_d$ \emph{kings graph} 
whose vertex set is the Cartesian product $[n_1] \times \cdots 
\times [n_d]$, where $[n]$ denotes $\{1,2,\ldots, n\}$. Two vertices 
are called \emph{neighbors} (or said to be \emph{adjacent}) when 
we can get from one to the other by a single generalized king's 
move. In other words, distinct vertices $\vx=(v_1,\ldots,v_d)$ and 
$\ux=(u_1,\ldots,u_d)$ in $K[n_1,\ldots,n_d]$ are neighbors if and only if 
$|v_i-u_i| \leq 1$ for each $i\in[d]$.  

The \emph{toroidal kings graph} $T[n_1,\ldots,n_d]$ is defined on 
the same vertex set as $K[n_1,\ldots,n_d]$, except that distinct 
vertices $\vx$ and $\ux$ are considered neighbors in 
$T[n_1,\ldots,n_d]$ if and only if $v_i-u_i \equiv -1$, $0$ or 
$1$ (mod $n_i$) for each $i\in [d]$.  The analysis of 
$T[n_1,\ldots,n_d]$ is much simpler than that of 
$K[n_1,\ldots\/,n_d]$ because vertices of $T[n_1,\ldots,n_d]$ 
have equally many neighbors. Note, however, that toroidal 
chessboards for which some $n_i<3$ sometimes require separate 
handling. In particular, the effects of adding 1 or $-1$ in 
coordinate $i$ are precisely the same when $n_i=2$.  For example, 
in $T[2,8,2]$ there are four ways to move from vertex $(2,6,1)$ 
to $(1,7,2)$: simply add any of the vectors $(1,1,1)$, 
$(-1,1,1)$, $(1,1,-1)$, $(-1,1,-1)$ to $(2,6,1)$. Thus $T[2,8,2]$ 
is really a \emph{multigraph}, in the sense that these two 
vertices are ``neighbors of multiplicity $4$.''%
\footnote{When $n_i$ is unity, each vertex is a multiple neighbor 
of itself. Removing index $i$ leads to an equivalent problem in 
lower dimensions, so we assume $n_i>1$ in this article.}
In $T[n_1,...,n_d]$ this \emph{multiplicity} is $2^c$, where $c$ 
is the number of coordinates $i$ at which two neighbors differ 
and for which $n_i=2$.  By counting multiplicities, all of our 
results can easily be extended to cover this situation, so we 
give it no further special treatment.

Now let $G$ be a general (loopless) graph with vertex set 
$V(G)$.  For a vertex $\vx\in V(G)$, $N(\vx)$ denotes the set of 
vertices adjacent to $\vx$. We call $N(\vx)$ the 
\emph{neighborhood} of $\vx$ in $G$, noting that $N(\vx)$ does 
not include the vertex $\vx$ itself.  Next, consider a whole 
number $k$ and a set $S\subseteq V(G)$ of vertices.  As 
introduced by Fink and Jacobson \cite{fj}, we say that $S$ is 
\emph{$k$-dependent} in $G$ if $|N(\vx) \cap S| \leq k$ for each 
$\vx\in S$, so that each vertex of $S$ has at most $k$ neighbors 
in $S$.  The name ``$k$-dependent'' arises from the case $k=0$, 
since a $0$-dependent set is usually called an \emph{independent} 
set in graph theory. The \emph{$k$-dependence number} of $G$, 
denoted by $\beta_k(G)$, is the maximum cardinality of a 
$k$-dependent set in $G$.

\begin{figure}[h]
\[
\begin{array}[c]{c}
  \underset{\mbox{\ \ (a) $k=4$, toroidal}}
    {\epsfig{file=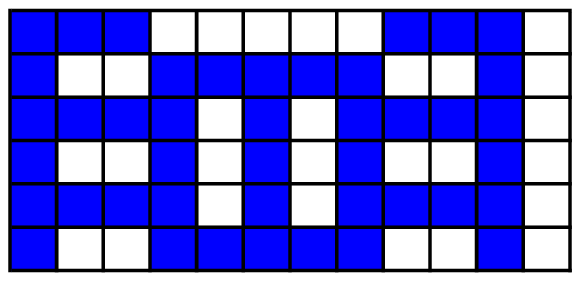,height=0.8in,width=1.5in}}\\ \\ \\ \\
  \underset{\mbox{(b) $k=5$, toroidal}}
    {\epsfig{file=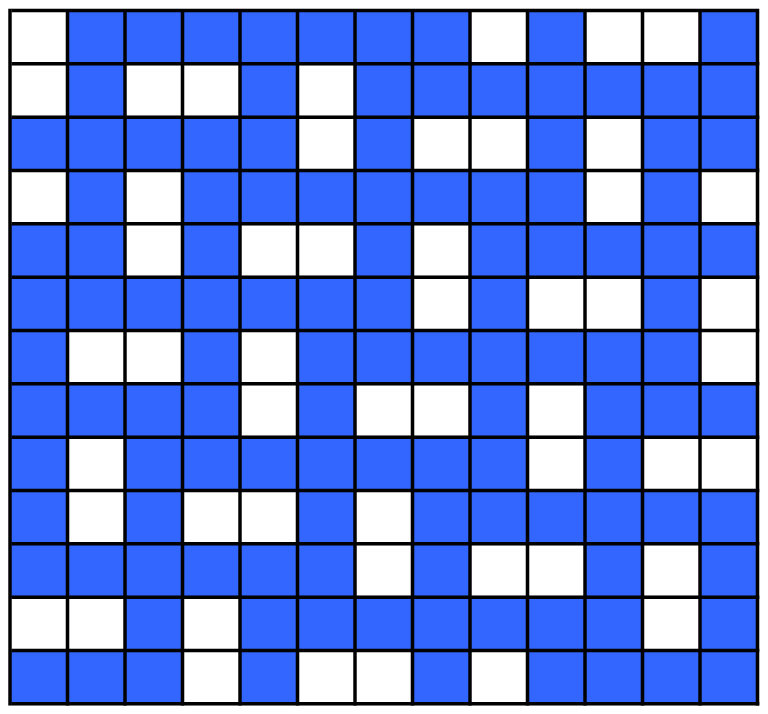,height=1.5in,width=1.5in}}
\end{array}
\qquad
\begin{array}[c]{c}
  \underset{\mbox{(c) half-dependent, non-toroidal}}
    {\epsfig{file=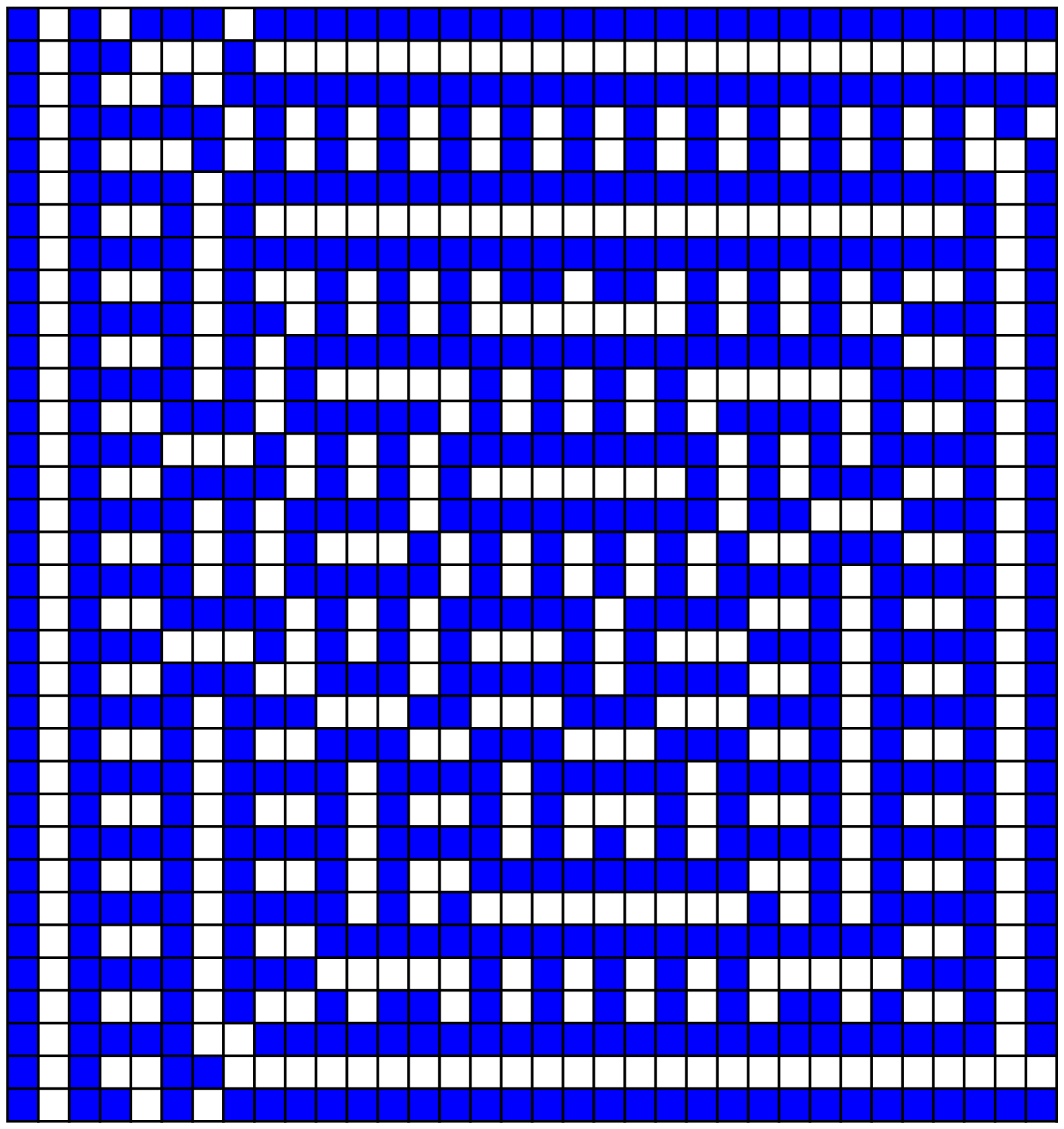,height=2.5in,width=2.5in}}\\ \\
  \underset{\mbox{\ \ \ \ \ (d) $k=4$, toroidal}}
    {\epsfig{file=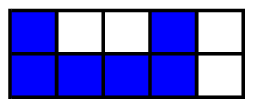,height=.5in,width=1in}}
\end{array}
\] 
\caption{Examples of $k$-dependent and half-dependent sets} 
\label{figure1}
\end{figure}

For a $k$-dependent set $S$ in a kings graph (toroidal or 
otherwise), we regard $S$ as the set of vertices or squares 
occupied by kings, no king having more than $k$ neighboring 
kings.  For example, Figure \ref{figure1}a shows a 4-dependent 
set of $43$ kings (indicated by dark squares) arranged in 
$T[6,12]$, proving that $\beta_4(T[6,12]) \geq 43$.  Likewise, 
Figure \ref{figure1}b shows a 5-dependent set of $117$ kings in 
$T[13,13]$, demonstrating that $\beta_5(T[13,13])\geq117$.

Our paper is partly motivated by \cite{dhlm}, which includes a 
section on ``$1/2$-domination'' of kings graphs $K[m,n]$ for 
small values of $m$. A subset $R$ of $V(G)$ is called a 
\emph{$1/2$-dominating set} if each vertex $\vx$ of $S = V(G) 
\setminus R$ satisfies $|R \cap N(\vx)| \geq |N(\vx)|/2$; in other 
words, $R$ is a \emph{dominating set} with the additional feature 
that each vertex not in $R$ is dominated by at least half of its 
neighbors. The \emph{$1/2$-domination number}, $\gamma_{1/2}(G)$, 
is the minimum cardinality among $1/2$-dominating sets in $G$.  
With our emphasis on $k$-dependence, we take the complementary 
perspective, defining a subset $S$ of $V(G)$ as 
\emph{half-dependent} in $G$  if each vertex $\vx$ of $S$ satisfies 
$|S \cap N(\vx)| \leq |N(\vx)|/2$.  The \emph{half-dependence 
number}, denoted by $h(G)$, is the maximum cardinality among 
half-dependent sets in $G$, so in general $h(G) = 
|V(G)|-\gamma_{1/2}(G)$.  The dark squares in Figure 
\ref{figure1}c form a half-dependent set of 694 vertices in 
$K[34,34]$, the white squares a $1/2$-dominating set, from which 
we have that $h(K[34,34]) \geq 694$ and $\gamma_{1/2}(K[34,34]) 
\leq 462$.  

For a graph $G$, we let $\tau_k(G)$ denote $\beta_k(G)/|V(G)|$, 
the maximum \emph{density} among $k$-dependent sets in $G$.  
Similarly, $\rho(G)=h(G)/|V(G)|$ denotes the maximum 
\emph{density} among half-dependent sets in $G$. For a given 
dimension $d$, $K^{(d)}[n]$ and $T^{(d)}[n]$ denote the special 
cases of $K[n_1,\ldots,n_d]$ and $T[n_1,\ldots,n_d]$, respectively, in 
which $n_i=n$ for each $i$. As a further shorthand, we use the 
following:
\begin{gather*}
  \tau^{(d)}_k(n)=\tau_k(T^{(d)}[n]),
    \quad \tau^{(d)}_k=\lim_{n\rightarrow\infty}\tau^{(d)}_k(n),\\
  \rho^{(d)}(n)=\rho(K^{(d)}[n]),
    \quad \rho^{(d)}=\lim_{n\rightarrow\infty}\rho^{(d)}(n).
\end{gather*}
In this paper, we prove that the limits $\tau^{(d)}_k$ and 
$\rho^{(d)}$ exist and we seek their exact values.  We provide 
good upper and lower bounds in many cases, and we obtain exact 
values for $\tau^{(2)}_k$ when $k\neq 4,5$.  Based on our 
results, we suspect that $\tau^{(d)}_k$ is a rational number for 
any $d$ and $k$.

For other results on combinatorial chessboard 
problems see \cite{bhh}, \cite{ffpp}, \cite{fricke}, \cite{hhr}, 
\cite{wrm}; for more on $k$-dependence see \cite{f}, 
\cite{glld}; and for other similar problems see \cite{g}, \cite{msv}.


\section{Two-Dimensional Kings Graphs} \label{2dkingsgraph}

The original motivation for the present paper concerns 
the following 
conjecture about the maximum density of kings on a standard 
non-toroidal, two-dimensional board.

\begin{conjecture}\label{3/5conjecture} For the case $d=2$, the
half-dependent limiting density is $\rho^{(2)}=3/5$ and the actual
half-dependent densities satisfy
\[ \rho^{(2)}(n) \geq \frac{3}{5}-\frac{C}{n^2}, \]
for all $n$ and some constant $C$.
\end{conjecture}

The example of a half-dependent set of 694 kings given in Figure 
\ref{figure1}c demonstrates 
\[ \rho^{(2)}(n) \geq \frac{3}{5}-\frac{2}{n^2}, \mbox{ when } n=34.\]
Indeed, we have computationally verified that $\rho^{(2)}(n) \geq 3/5-7/(5n^2)$ 
for each $n \leq 35$.  In a later section, we give more credence 
to Conjecture \ref{3/5conjecture} by establishing the upper bound 
$\rho^{(2)}\leq0.608956$.  Meanwhile, we can give some lower 
bounds on $\rho^{(2)}(n)\cdot n^2=h(K[n,n])$.

\begin{theorem}\label{43/72exampleUsed}
The maximum size of a half-dependent set in the $n \times n$ kings 
graph satisfies the following lower bounds, for some constant $C$:
\begin{itemize}
\item
$h(K[n,n]) \geq {\displaystyle \frac{3}{5}n^2-\frac{n}{30}-C}$,
          if $n\equiv 0 \pmod5$;
\item
$h(K[n,n]) \geq {\displaystyle \frac{3}{5}n^2-\frac{4n}{30}-C}$,
          if $n\equiv 1 \pmod5$;
\item
$h(K[n,n]) \geq {\displaystyle \frac{3}{5}n^2-\frac{2n}{30}-C}$,
          if $n\equiv 2 \pmod5$;
\item
$h(K[n,n]) \geq {\displaystyle \frac{3}{5}n^2+\frac{3}{5}}$,
          \qquad\; if $n\equiv 3 \pmod5$;
\item
$h(K[n,n]) \geq {\displaystyle \frac{3}{5}n^2-\frac{3n}{30}-C}$,
          if $n\equiv 4 \pmod5$.
\end{itemize}
\end{theorem}

\begin{figure}[t]
\[\underset
    {\parbox{.57in}{\footnotesize(a) $n=3$\\ \phantom{\ \ }6 kings}}
    {\epsfig{file=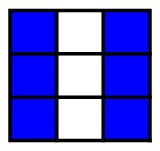,height=.55in,width=.55in}} \ \
  \underset
    {\parbox{.6in}{\footnotesize(b) $n=4$\\ \phantom{\ \ }9 kings}}
    {\epsfig{file=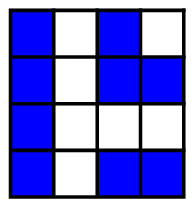,height=0.7in,width=0.7in,angle=90}} \ \
  \underset
    {\parbox{.7in}{\footnotesize(c) $n=5$\\ \phantom{\ \ }15 kings}}
    {\epsfig{file=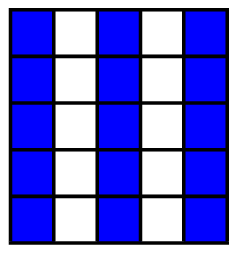,height=.9in,width=.9in,angle=90}} \ \
  \underset
    {\parbox{.7in}{\footnotesize(d) $n=6$\\ \phantom{\ \ }22 kings}}
    {\epsfig{file=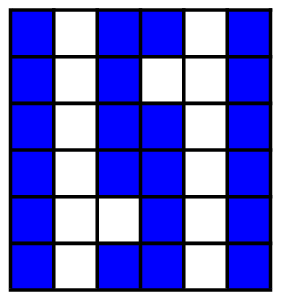,height=1.05in,width=1.05in,angle=90}} \ \
  \underset
    {\parbox{.7in}{\footnotesize(e) $n=7$\\ \phantom{\ \ }28 kings}}
    {\epsfig{file=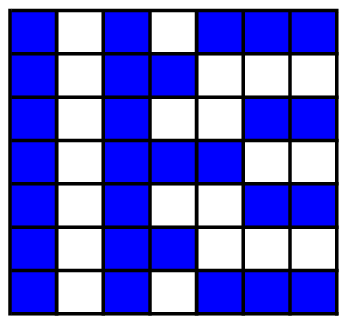,height=1.2in,width=1.2in}}
\]
\[\underset
    {\parbox{.7in}{\footnotesize(f) $n=8$\\ \phantom{\ \ }39 kings}}
    {\epsfig{file=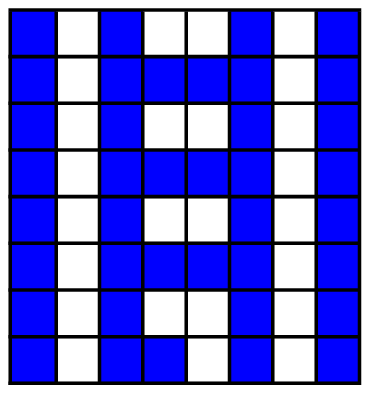,height=1in,width=1in,angle=-90,origin=c}}\ \
  \underset
    {\parbox{.7in}{\footnotesize(g) $n=9$\\ \phantom{\ \ }49 kings}}
    {\epsfig{file=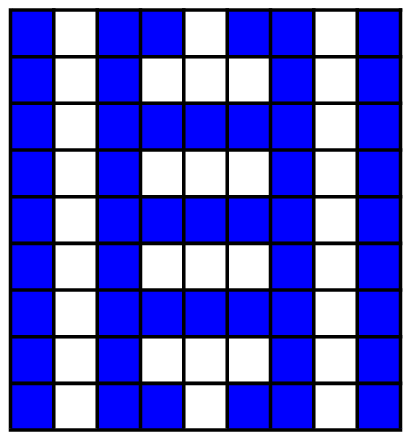,height=1.1in,width=1.1in,angle=90}}\ \
  \underset
    {\parbox{.7in}{\footnotesize(h) $n=10$\\ \phantom{\ \ }59 kings}}
    {\epsfig{file=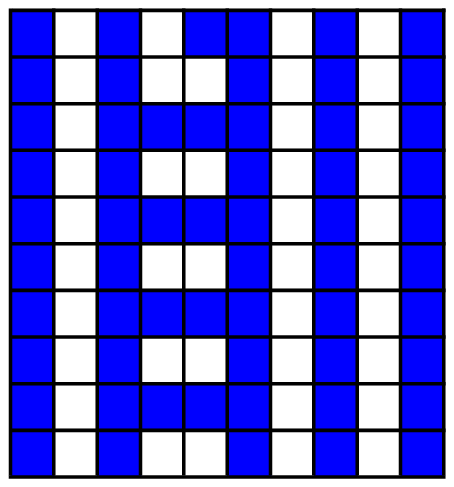,height=1.25in,width=1.25in,angle=90}}\ \
  \underset
    {\parbox{.7in}{\footnotesize(i) $n=11$\\ \phantom{\ \ }73 kings}}
    {\epsfig{file=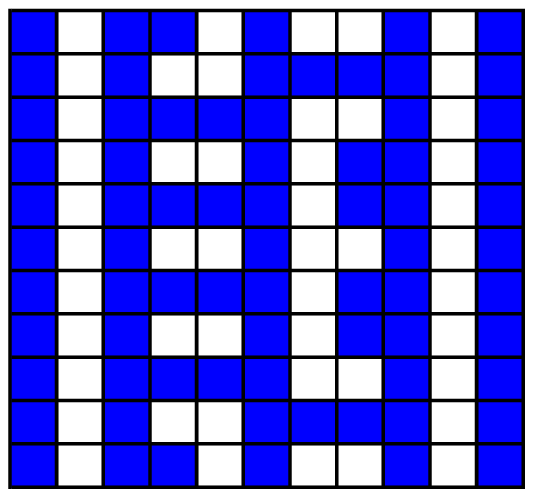,height=1.35in,width=1.35in}}
\]
\caption{Examples of maximum-density, half-dependent sets.} 
\label{halfdependent}
\end{figure}

\begin{proof}
Let $\mathcal{C}$ and $\mathcal{D}$ denote the toroidal 
arrangements in Figures \ref{figure1}d and \ref{figure1}a, 
respectively.  Upon stacking $m$ copies of 
$\mathcal{C}$ and removing the king at $(2m,3)$,%
\footnote{The reader should note that all illustrations in this 
paper follow matrix indexing, so that vertex $(i,j)$ appears as a 
square in row $i$, column $j$.}
we obtain an arrangement $\mathcal{A}'$ comprised of $6m-1$ kings 
in $K[2m,5]$. In $K[2m+1,5]$ upon placing a copy of 
$\mathcal{A}'$ in $\{2,3,\ldots,2m+1\} \times [5]$ and including 
additional kings at $(1,1), (1,2)$ and $(1,4)$, we have an 
arrangement $\mathcal{A}''$ comprised of $6m+2$ kings.  Thus via 
$\mathcal{A}'$ or $\mathcal{A}''$ we have an arrangement 
$\mathcal{A}$ in $K[n,5]$ comprised of $3n-1$ kings, using no 
kings in column 5.

Similarly, stacking $m= \lfloor(n-2)/6\rfloor$ copies of 
$\mathcal{D}$ and placing the result within $\{2,3,\ldots,6m+1\} 
\times [12]$ we obtain an arrangement $\mathcal{B}$ in $K[n,12]$ 
comprised of $43m = (43/72)(12n) -c=(43n/6)-c$ kings for some 
constant $c$ (based on the fact that we have generously left row 
1 and rows $6m+2$ through $n$ devoid of kings). Arrangement 
$\mathcal{B}$ has no kings in column 12.

We now use $\mathcal{A}$ and $\mathcal{B}$ to construct the 
desired half-dependent arrangements within $K[n,n]$.  If $n=5m+3$ 
then place $m$ copies of $\mathcal{A}$ side by side in $[n] 
\times \{3,4,\ldots,n-1\}$ and kings everywhere in columns 1 and $n$ 
to verify the result, as in Figure \ref{halfdependent}f when 
$n=8$. If $n=5m+15$ then place $m$ copies of $\mathcal{A}$ side 
by side in $[n] \times \{3,4,\ldots,n-13\}$ followed by a copy of 
$\mathcal{B}$ in columns $n-12$ through $n-1$ and kings 
everywhere in columns 1 and $n$, and similarly using two copies 
of $\mathcal{B}$ if $n=5m+27$, three copies of $\mathcal{B}$ if 
$n=5m+39$, and four copies of $\mathcal{B}$ if $n=5m+51$.  Thus 
in each congruence case the result is verified for large $n$, and 
small values of $n$ are automatically correct by specifying the 
constant $C$ sufficiently large in compensation.
\end{proof}

\section{Limiting Densities}

Our next result shows that the limiting densities exist and 
illustrates the tight relationship between the half-dependent 
non-toroidal and $k$-dependent toroidal problems. It also 
provides a lower bound on the limiting densities.

\begin{theorem}\label{existence}
The limiting values $\tau^{(d)}_k$ and $\rho^{(d)}$ exist, and 
satisfy
\begin{enumerate}
\item 
  $\tau^{(d)}_k\geq\tau_k(T[n_1,\ldots,n_d])$, for any choice of $d$, 
  $k$, and $n_1,\ldots,n_d$;
\item
  $\tau^{(d)}_k=\rho^{(d)}$, in the particular case where
  $k=(3^d-1)/2$.
\end{enumerate}
\end{theorem}

\begin{figure}
\[
  \epsfig{file=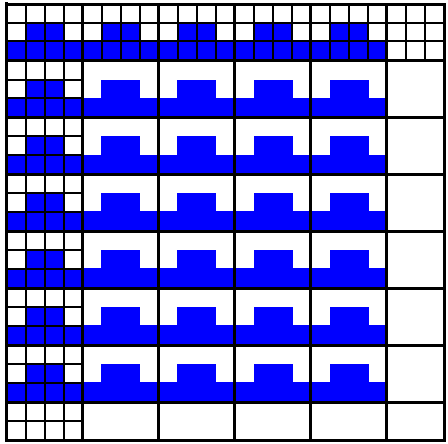,height=2in,width=2in}
\]
\caption{ A packing of $T[23,23]$ with $\lfloor 
23/3\rfloor\cdot\lfloor 23/4\rfloor$ copies of $T[3,4]$, each 
copy containing a $4$-dependent set of maximum density. } 
\label{packing}
\end{figure}

\begin{proof}
Consider $n_1,\ldots,n_d>0$ and $n>\max_in_i$.  The 
Quotient-Remainder Theorem (division algorithm) allows us to 
uniquely write $n=n_i\lfloor n/n_i\rfloor +r_i$ for some 
$r_i\in\{0,\ldots,n_i-1\}$.  We can then pack $T^{(d)}[n]$ with 
$\prod_{i=1}^d \lfloor n/n_i\rfloor$ non-overlapping copies of 
$T[n_1,\ldots,n_d]$. These copies can be aligned so that the 
toroidal boundaries are compatible from one copy to the next, 
except for those abutting the ``remainder'' sections of length 
$r_i$ in each coordinate. An example of such a packing is 
illustrated in Figure \ref{packing}.  Next, we place a 
$k$-dependent set of density $\tau_k(T[n_1,\ldots,n_d])$ within 
each copy of $T[n_1,\ldots,n_d]$. This yields a $k$-dependent set 
$S\subseteq T^{(d)}[n]$, thereby giving us the bound
\begin{eqnarray}
\tau^{(d)}_k(n)
  \geq \frac{|S|}{n^d}
  &=& \frac
        {\beta_k(T[n_1,\ldots,n_d])\prod_{i=1}^d\lfloor n/n_i\rfloor}
        {n^d} \nonumber\\
  &=& \tau_k(T[n_1,\ldots,n_d])
       \prod_{i=1}^d \frac{\lfloor n/n_i\rfloor}{n/n_i}.\label{taucopy}
\end{eqnarray}
In the special case where $n_i=m<n$ for all $i$, this implies
\[ \tau_k^{(d)}(n)
   \geq \tau_k^{(d)}(m)\left(\frac{\lfloor n/m\rfloor}{n/m}\right)^d.
\]
Taking the limit infimum as $n\to\infty$ yields
\begin{equation}\label{tauliminf} 
  \liminf_{n\to\infty}\tau_k^{(d)}(n)
  \geq \tau_k^{(d)}(m),\quad\forall m>0.
\end{equation}
From here, the limit supremum as $m\to\infty$ gives us
\[ \liminf_{n\to\infty}\tau_k^{(d)}(n)
   \geq \limsup_{m\to\infty}\tau_k^{(d)}(m).
\]
Consequently, $\tau_k^{(d)}=\lim_{n\to\infty}\tau_k^{(d)}(n)$ exists.
Combined with inequalities (\ref{taucopy}) and (\ref{tauliminf}), 
this also proves statement (a).

Now let $k=(3^d-1)/2$.  By deleting the ``boundary'' kings from a 
$k$-dependent subset of density $\tau_k^{(d)}(n)$ on the toroidal 
board $T^{(d)}[n]$, we obtain a half-dependent subset $S$ of the 
non-toroidal board $K^{(d)}[n]$.  This implies that
\[ \rho^{(d)}(n)
   \geq \frac{|S|}{n^d}
   \geq \tau_k^{(d)}(n) - \frac{n^d-(n-2)^d}{n^d}.
\]
Reversing the roles of the two boards yields the analogous 
inequality
\[ \tau_k^{(d)}(n) \geq \rho^{(d)}(n) - \frac{n^d-(n-2)^d}{n^d}. \]
Combining these and taking the limit proves statement (b) and the 
existence of $\rho^{(d)}$.
\end{proof}

The analogue of Theorem \ref{existence}(a) in which 
$\tau^{(d)}_k$ and $\tau_k(T[n_1,\ldots,n_d])$ are replaced by 
$\rho^{(d)}$ and $\rho^{(d)}(K[n_1,\ldots,n_d])$ fails: Theorems
\ref{43/72exampleUsed} and \ref{maintheorem} show that 
$\rho^{(2)}\leq0.609<\rho^{(2)}(3)=2/3$. 

Figure \ref{figure1}d shows a 4-dependent set of 6 kings in 
$T[2,5]$, where each king has four neighboring kings when the 
neighbors are counted with multiplicity.  Combining Theorem 
\ref{existence} with the examples of Figures \ref{figure1}b and 
\ref{figure1}d, we have verified the following numerical lower 
bounds.

\begin{corollary}\label{lowerbounds}
For 4- and 5-dependent kings graphs in two dimensions, we have 
the lower bounds $\rho^{(2)} = \tau^{(2)}_4 \geq 3/5$
and $\tau^{(2)}_5 \geq 9/13$.
\end{corollary}


\section{Binary Linear Programming}\label{binaryLP}

The $k$-dependence number $\beta_k(G)$ of a given graph $G$ can 
be computed, in principle, by reformulating the corresponding 
maximization problem.  To each $k$-dependent set $S$, we associate 
the \emph{characteristic vector} with entries
\[ x_\vx =\begin{cases}
   1, & \text{if $\vx\in S$,}\\
   0, & \text{otherwise.}
   \end{cases}
\]
To each vertex $\vx\in V(G)$, we associate the inequality%
\footnote{In the case of a multigraph, such as a toroidal kings 
graph for which with some $n_i=2$, the summation term in 
(\ref{constraint}) must be modified to account for 
multiplicities.}
\begin{equation}
  \left(|N(\vx)|-k\right)x_\vx+\sum_{\ux\in N(\vx)}x_\ux \leq |N(\vx)|
\label{constraint}
\end{equation}
and impose the restriction that $x_\vx$ be $0$ or $1$. Separate 
consideration of the cases $x_\vx=0$ and $x_\vx=1$ shows that the 
resulting system of inequality constraints precisely describes 
$k$-dependent sets.  The optimization problem consists of 
maximizing the linear function $\sum_{\vx\in V(G)}x_\vx$ subject 
to this system of inequalities and the $0-1$ restrictions.  This 
is an example of a \emph{binary linear programming} (binary LP) 
problem. The optimal value in this problem is $\beta_k(G)$ and an 
optimal solution vector corresponds to a maximum $k$-dependent 
set.

The optimization problem for determining $h(G)$ can be formulated 
similarly, except that the constraint associated with each vertex 
$\vx$ becomes
\[ \left\lceil \frac{|N(\vx)|}{2} \right\rceil x_\vx
   + \sum_{\ux\in N(\vx)}x_\ux \leq |N(\vx)|.
\]
Applying a binary LP solver to this problem, we determined the 
values of $h(K[n,n])$ shown in Table \ref{halftable} for $1 
\leq n \leq 11$, along with the sample half-dependent sets of 
$h(K[n,n])$ kings shown in Figure \ref{halfdependent}.  For 
$n=8,9,10,11,$ all optimal patterns look like the sample shown; 
for $n=7$, there are several distinct optimal patterns, including 
a pattern consisting of vertical stripes. 

\begin{table}
\footnotesize\[
\begin{tabular}{|c|c|c|c|c|c|c|c|c|c|c|c|} \hline
$n $         & 1 & 2& 3 &4 &5&6&7&8&9&10&11 \\
\hline
$h(K[n,n])$   & 1 & 2& 6 & 9& 15&22& 28&39&49&59&73\\
\hline
\end{tabular}
\]
\caption{Maximum number of kings in half-dependent sets.} 
\label{halftable}
\end{table}

The ``binary'' restriction of these LPs (namely, that $x_\vx$ be 
$0$ or $1$) implies the possibility of searching through all 
$2^{|V(G)|}$ arrangements, looking for a largest $k$-dependent 
one.  Although modern software for solving LPs manages to avoid 
considering nearly so many arrangements, there quickly comes a 
point where $|V(G)|$ is simply too large for the above method to 
be practical.  The next section shows how the binary LP 
perspective can still allow for efficient calculation of good 
upper bounds.


\section{Upper Bounds and Exact Values}\label{upperbound}

In this section, we calculate upper bounds by solving binary LPs 
on relatively small vertex sets.  In fact, the system of 
inequalities derived in the preceding section can lead to general 
upper bounds without even having to solve the associated linear 
program! As an example, consider the problem of calculating 
$\beta_2(T[n,n])$, the maximum number of kings that can be placed 
on an $n \times n$ toroidal board with no king having more than 
two neighboring kings.  In this case, the LP constraint 
(\ref{constraint}) is
\[ 6x_\vx+\sum_{\ux\in N(\vx)}x_\ux \leq 8. \]
Note that $x_\vx$ appears in nine of these $n^2$ constraints: once 
(with coefficient 6) in its own associated constraint and once 
(with coefficient 1) in each of the eight constraints associated 
with its neighbors. Summing the constraints over $\vx$, we find 
that the characteristic vector of a 2-dependent set $S$ satisfies
\[ 14|S| = 14\sum_{\vx\in V(G)}x_\vx \leq 8n^2.\]
Thus $|S| \leq (4/7)n^2$, which implies that $\tau^{(2)}_2(n) 
\leq 4/7$.

But we can do better. Observe that, for any 2-dependent set $S$ 
and any vertex $\vx$ in an $n \times n$ toroidal board, we have 
$\sum_{\ux\in N(\vx)}x_\ux \leq 6$, since placing kings at 7 of the 8 
neighbors of $\vx$ always violates 2-dependence.  This allows us to 
re-derive the LP constraint associated with $\vx$, so that it is
replaced by 
\[ 4x_\vx+\sum_{\ux\in N(\vx)}x_\ux \leq 6.\]
Summing this new set of constraints yields $12|S| \leq 6n^2$, 
thereby improving the bound to $\tau^{(2)}_2(n) \leq 1/2$.  On 
the other hand, we can form a 2-dependent set by placing a king 
at $\vx=(v_1,v_2)$ if and only if $v_2$ is even.  This shows that 
$\beta_2(T[n,n])\geq n\lfloor n/2 \rfloor$ and therefore $\lfloor 
n/2\rfloor/n \leq \tau^{(2)}_2(n) \leq 1/2$.  In the limit, we 
obtain $\tau^{(2)}_2 = 1/2$. There are simpler ways to obtain the 
exact value for $\tau^{(2)}_2$, but the approach just given can 
be generalized, as we show next.

An \emph{automorphism} of a graph $G$ is a bijection 
$f:V(G)\rightarrow V(G)$ that preserves adjacency, so that 
neighbors are mapped to neighbors.  Graph $G$ is called 
\emph{vertex-transitive} if for every two vertices $\vx,\ux$ 
there exists an automorphism $f$ for which $f(\vx)=\ux$.  In 
other words, $G$ is vertex-transitive when each vertex plays the 
same structural role in $G$ as any other vertex, such as happens 
in toroidal kings graphs but not in kings graphs.  Note that in a 
vertex-transitive graph, the neighborhoods $N(\vx)$ all have the 
same cardinality.  In the following result, $\langle V' \rangle$ 
denotes the subgraph \emph{induced} by a subset $V'$ of $V(G)$, 
namely, the subgraph of $G$ formed by deleting all vertices of 
$G$ not in $V'$.  

\begin{proposition}\label{nbhdprop}
Suppose $G$ is a vertex-transitive graph and let $\beta^*_k$ 
denote the quantity $\beta_k(\langle N(\vx) \rangle)$, which is 
independent of the choice of vertex $\vx$. Then
\[ \beta_k(G) \leq \frac{\beta^*_k |V(G)|}{\beta^*_k -k+|N(\vx)|}. \]
\end{proposition}

\begin{proof}
Associate a constraint
$(\beta^*_k-k)x_\vx+\sum_{\ux\in N(\vx)}x_\ux \leq \beta^*_k$ to each vertex
$\vx\in V(G)$, sum the constraints over all $\vx$, and deduce the maximum
of $|S|=\sum_{\vx\in V(G)}x_\vx$.
\end{proof}

As an application of Proposition \ref{nbhdprop}, we return to the 
two-dimensional setting of arranging kings on a large $n \times 
n$ board.  In this case, the values of $\beta^*_k$ for 
$G=T^{(2)}[n]$ are easily calculated by hand:
\[
   \beta^*_0=\beta^*_1=4,\qquad
   \beta^*_2=\beta^*_3=6,\qquad
   \beta^*_4=\beta^*_5=\beta^*_6=\beta^*_7=\beta^*_8=8.
\]
The corresponding upper bounds on $\tau^{(2)}_k$ are
\begin{gather*}
  \tau_0^{(2)} \leq \frac{4}{12},\quad
  \tau_1^{(2)} \leq \frac{4}{11},\quad
  \tau_2^{(2)} \leq \frac{6}{12},\quad
  \tau_3^{(2)} \leq \frac{6}{11},\quad
  \tau_4^{(2)} \leq \frac{8}{12},\\
  \tau_5^{(2)} \leq \frac{8}{11},\quad
  \tau_6^{(2)} \leq \frac{8}{10},\quad
  \tau_7^{(2)} \leq \frac{8}{9},\quad
  \tau_8^{(2)} \leq \frac{8}{8}\,,
\end{gather*}
whereas the known exact values for $\tau^{(2)}_k$ turn out to be
\[ 
  \tau^{(2)}_0=\frac{1}{4}\,,\quad 
  \tau^{(2)}_1=\frac{1}{3}\,, \quad
  \tau^{(2)}_2=\tau^{(2)}_3=\frac{1}{2}\,,\quad
  \tau^{(2)}_6=\frac{4}{5}\,, \quad
  \tau^{(2)}_7=\frac{8}{9}\,,\quad
  \tau^{(2)}_8=1\,.
\]
Thus, the upper bound of Proposition \ref{nbhdprop} is tight for 
$k\in\{0,2,6,7,8\}$.  We can verify this tightness on a 
case-by-case basis.  Note that $\tau^{(2)}_2=1/3$ was proved 
earlier in this section.  The values of $\tau^{(2)}_k$ with 
$k\in\{0,7,8\}$ correspond to taking $d=2$ in the more general 
formulas
\[ \tau^{(d)}_0=2^{-d}, \quad \tau^{(d)}_{3^d-2}=1-3^{-d}, \quad
   \tau^{(d)}_{3^d-1}=1,
\]
which can be established easily.  The situation for $k=6$ also 
admits a nice generalization to higher dimensions, but is 
somewhat trickier; see section \ref{congruences}.    

The upper bound in Proposition \ref{nbhdprop} is not tight for 
$k\in\{1,3,4,5\}$, so additional methods are needed. In section 
\ref{k=4&5}, we generalize Proposition \ref{nbhdprop} in a way 
that improves the upper bound in these four cases.  The exact 
values for $\tau^{(2)}_1$ and $\tau^{(2)}_3$ are obtained by 
other means in section \ref{others}.


\section{Linear Congruences and Lower Bounds}\label{congruences}

Here we use linear congruences to produce specific $k$-dependent 
sets in $T^{(d)}[n]$, thus giving lower bounds for $\tau_k^{(d)}$.
These will match the upper bounds derived in the preceding 
section.  For this purpose, we regard the coordinates of vertices 
in $T^{(d)}[n]$ as elements of the group $\Z_n$ of integers 
modulo $n$. For a vector $\cx\in (\Z_n)^d$ and a set $R\subseteq 
\Z_n$, let $S(n,d,\cx,R)$ denote the vertex set 
\[ S(n,d,\cx,R)=\{\vx\in V(T^{(d)}[n]) :
   \cx\cdot \vx\in R\}.
\]
For which values of $k$ is $S(n,d,\cx,R)$ a $k$-dependent 
set? Consider a vertex $\vx\in S(n,d,\cx,R)$
and a nontrivial vector $\yx\in \{-1,0,1\}^d$.  Note that
the neighbor $\vx+\yx$ of $\vx$ belongs to
$S(n,d,\cx,R)$ if and only if
$\cx\cdot \vx+\cx\cdot \yx \in R$.
Thus, if we let $f(\cx,r)$ denote the number of 
nontrivial vectors $\yx\in \{-1,0,1\}^d$ for which 
$\cx\cdot \yx+r\in R$, then
$f(\cx,\cx\cdot \vx)$ is the number of elements of $S(n,d,\cx,R)$
neighboring $\vx$.  We therefore have the following:

\begin{lemma}
For the choice $k =\max_{r \in R}f(\cx,r)$, the set 
$S(n,d,\cx,R)$ is $k$-dependent in $T^{(d)}[n]$.
\end{lemma}

To obtain exact values of maximum toroidal densities in certain 
cases, we can combine the above lemma with the following 
well-known and easily verified fact:  to each integer $m \neq 0$ 
there corresponds exactly one positive integer $d$ and one vector 
$\yx \in \{-1,0,1\}^d$ for which $m =\cx\cdot \yx$, where $\cx = 
(3^{d-1},3^{d-2},\ldots,3^1,3^0)$.  For instance, corresponding 
to $m=19$ are the choices $d=4$ and $\cx=(1,-1,0,1)$, based on 
the fact that 19 is expressible as $27-9+1$ using sums and/or 
differences
of distinct powers of 3.%
\footnote{The representation of $m$ as $\cx\cdot \yx$
uses \emph{balanced ternary notation}, which is much like base 3,
or ternary, notation.  To derive the representation given here for
$m=19$, rewrite the standard ternary form as 
$19 = (2)3^2+(0)3^1+(1)3^0
    = (3-1)3^2+(0)3^1+(1)3^0
    =(1)3^3+(-1)3^2+(0)3^1+(1)3^0
$.}  This leads to the following result on densities.

\begin{theorem}\label{3d-3}
For $k=3^d-3$, we have
\[ \tau_k^{(d)} = \frac{3^d-1}{3^d+1}. \]
In particular, $\tau^{(2)}_6 =4/5$ and $\tau_{24}^{(3)} = 13/14$.
\end{theorem}

\begin{proof}
Consider $S = S(n,d,\cx,R)$ with $n = (3^d+1)/2$, 
$R = \Z_n \setminus \{0\}$, and
\[ \cx = (3^{d-1},3^{d-2},\ldots,3^1,3^0), \]
as in the above lemma. It is straightforward to verify that, for 
any $r\in R$, there exist nontrivial vectors
$\yx,\yx' \in \{-1,0,1\}^d$ for which $\cx\cdot \yx=-r$ 
and $\cx\cdot \yx'=n-r$. Thus $f(\cx,r) 
\leq 3^d-3$, since at least two of the $3^d-1$ choices of 
nontrivial vectors in $\{-1,0,1\}^d$ must have
$\cx\cdot \yx+r \notin R$. Therefore, $S$ is $(3^d-3)$-dependent.

Now observe that, for each choice of values 
$v_1,v_2,\ldots,v_{d-1}$, there is exactly one choice of $v_d \in 
\Z_n$ for which $\vx=(v_1,v_2,\ldots,v_d)\notin S$.  
Consequently,
\[ |S| = \frac{|R|}{n} = \frac{3^d-1}{3^d+1}, \]
so Theorem \ref{existence}(a) gives us 
\[ \tau_{k}^{(d)} \geq \tau_{k}(T^{(d)}[n])
   \geq \frac{3^d-1}{3^d+1}.
\]
For the matching upper bound, we apply Proposition \ref{nbhdprop}
with $G=T^{(d)}[n]$, $k=3^d-3$, and $|N(v)|=3^d-1$ to get
\begin{eqnarray*}
\tau_{k}^{(d)}(n)
  = \frac{\beta_k(G)}{|V(G)|}
  &\leq& \frac{\beta_k^*}{\beta_k^*-k+|N(v)|}\\
  &\leq& \frac{3^d-1}{(3^d-1)-(3^d-3)+(3^d-1)}
  = \frac{3^d-1}{3^d+1},\\
\end{eqnarray*}
where we use the fact that $\beta_k^* \leq |N(v)|$. Therefore
$\tau_{k}^{(d)}\leq (3^d-1)/(3^d+1)$.
\end{proof}


\section{Improving Upper Bounds for $\tau^{(2)}_4$ and
$\tau^{(2)}_5$}\label{k=4&5}

We now generalize the binary LP approach of section 
\ref{upperbound} to address the unsolved problems of determining 
values for $\tau^{(2)}_4$ and $\tau^{(2)}_5$.  Recall that in 
Proposition \ref{nbhdprop} we used the inequality 
$(\beta^*_k-k)x_\vx+\sum_{\ux\in N(\vx)}x_\ux \leq \beta^*_k$, which is 
valid for any $k$-dependent set in a vertex-transitive graph.  
Here we seek other inequality constraints that are valid for all 
$k$-dependent sets.

Consider any \emph{weighting function} $\omega : V(G) \rightarrow
[0,\infty)$, not everywhere zero, and let $W(\omega)$ denote the 
\emph{total weight}, $\sum_{\vx\in V(G)}\omega(\vx)$, of $\omega$. 
For a given $k$ and $\omega$, let $M_k(G,\omega)$ denote the 
maximum value of $\sum_{\vx\in V(G)}\omega(\vx) \cdot x_\vx$ over all 
$k$-dependent sets $S$ in $V(G)$. To compute $M_k(G,\omega)$ we 
simply maximize the objective function $\sum_{\vx\in V(G)}\omega(\vx) 
\cdot x_\vx$ (instead of $\sum_{\vx\in V(G)}x_\vx$), using the same 
constraints as when computing $\beta_k(G)$ in section 
\ref{binaryLP}.

We have already seen two examples of such weighting functions. 
One is the case where $\omega$ is the constant function 
$\omega(\vx)=1$ (for all $\vx$), in which case $M_k(G,\omega)$ 
corresponds to the value of $\beta_k(G)$.  The other example is
\[  \omega(\vx)=
    \begin{cases}
       1,&\text{if $\vx \in N(\ux)$},\\
       \beta^*_k-k,&\text{if $\vx=\ux$,}\\
       0,&\text{otherwise},
    \end{cases}
\]
for some fixed vertex $\ux$; in this case $M_k(G,\omega)$ equals 
$\beta^*_k$.  In general, for any $\omega$ and any $k$-dependent
set in $V(G)$, we always have
\[ \sum_{\vx\in V(G)}\omega(\vx) \cdot x_\vx \leq M_k(G,\omega), \]
simply by the definition of $M_k(G,\omega)$.

\begin{lemma}\label{weightedtransitive}
Consider any weighting function $\omega$ on a vertex-transitive 
graph $G$.  Then an upper bound for the maximum density among 
$k$-dependent sets in $G$ is given by
$\tau_k(G) \leq M_k(G,\omega)/W(\omega)$.
\end{lemma}

\begin{proof}
Let $\Gamma$ denote the group of all automorphisms on $G$ and let 
$F$ denote $|\{f\in\Gamma:f(\vx)=\vx\}|$, a number which is 
independent of the choice of $\vx\in V(G)$.  Consider a weighting 
function $\omega$ and a $k$-dependent set $S$ in $V(G)$.  If $x_\vx$ 
is the characteristic vector for $S$ and $f$ is some 
automorphism, then the vector $x'_\vx=x_{f(\vx)}$ is also the 
characteristic vector of a $k$-dependent set.  Therefore,
\[ 
  \sum_{\vx\in V(G)}\omega(\vx) \cdot x_{f(\vx)} \leq M_k(G,\omega)
\]
holds for each automorphism $f\in\Gamma$.  Summing these 
inequalities, one per automorphism $f$, we see that each variable
$x_\vx$ appears on the left-hand side with total coefficient equal 
to $W(\omega)F$.  Thus we obtain the inequality
\[ W(\omega)F|S| \leq M_k(G,\omega) |\Gamma|. \]
Using the fact (see p. 89 of \cite{hun}, for example) that 
$|\Gamma|=F|V(G)|$, we obtain
\[ |S| \leq \frac{M_k(G,\omega)}{W(\omega)} |V(G)|, \]
proving the claim.
\end{proof}

The preceding lemma gives a very general tool for upper bounding 
the toroidal limiting density by means of relatively small 
non-toroidal kings graphs, as shown in the next result.  Note 
that this is the only place where we explicitly consider 
$k$-dependence on a \emph{non-toroidal} kings graph.

\begin{lemma}\label{weightedkings}
For any weighting function $\omega$ on $K[n_1,\ldots,n_d]$, we have
\[
  \tau_k^{(d)} \leq \frac{M_k(K[n_1,\ldots,n_d],\omega)}{W(\omega)}.
\]
Furthermore,
$\inf_\omega \left[M_k(K^{(d)}[n],\omega)/W(\omega)\right]
  \to\tau_k^{(d)}$ as $n\to\infty$.
\end{lemma}

\begin{proof}
Consider any weighting function $\omega$ for $K[n_1,\ldots,n_d]$.
Then for all $n > n_1,n_2,\ldots,n_d$, consider the weighting 
function $\omega'$ on $T^{(d)}[n]$ defined by
\[ \omega'(\vx) = 
   \begin{cases}
     \omega(\vx),& \text{if $\vx \in V(K[n_1,\ldots,n_d])$,}\\
     0,&\text{otherwise.}\\
   \end{cases}
\] 
The intersection of each $k$-dependent subset of $V(T^{(d)}[n])$ with 
$V(K[n_1,\ldots,n_d])$ is $k$-dependent in $K[n_1,\ldots,n_d]$, so 
$M_k(T^{(d)}[n],\omega') = M_k(K[n_1,\ldots,n_d],\omega)$. Because 
$T^{(d)}[n]$ is vertex-transitive, Lemma \ref{weightedtransitive} 
implies that
\[ \tau_k^{(d)}(n)
   \leq \frac{M_k(T^{(d)}[n],\omega')}{W(\omega')}
   = \frac{M_k(K[n_1,\ldots,n_d],\omega)}{W(\omega)}
\]
for all suitably large $n$.  Letting $n\rightarrow \infty$ 
completes the proof of the inequality.

To prove the second statement, consider the case where $n_i=n$ for all
$i$ and $\omega$ is the constant function $\omega(\vx)\equiv1$.  If $S$
is a $k$-dependent set of maximum density in $K^{(d)}[n]$, then its 
characteristic vector maximizes $\sum_{\vx}\omega(\vx)x_\vx$.  Viewing 
$K^{(d)}[n]$ as a subset of $T^{(d)}[n+1]$, we obtain
\[ \frac{|S|}{(n+1)^d}
   \leq \tau_k^{(d)}
   \leq \frac{M_k(K^{(d)}[n],\omega)}{W(\omega)}
   = \frac{|S|}{n^d}.
\]
This implies that
\[ \tau_k^{(d)}
   \leq \frac{M_k(K^{(d)}[n],\omega)}{W(\omega)}
   \leq \tau_k^{(d)} \left(\frac{n+1}{n}\right)^d,
\]
so that $M_k(K^{(d)}[n],\omega)/W(\omega)$ can be made 
arbitrarily close to $\tau_k^{(d)}$ by choosing $n$ sufficiently 
large.
\end{proof}

We are now ready for our main result.

\begin{theorem}\label{maintheorem} The limiting densities for 4- and
5-dependent toroidal kings graphs in two dimensions satisfy the 
bounds
\[ 0.6=\frac{3}{5} \leq \tau^{(2)}_4 
   \leq 0.608956 
\]
and
\[ 0.6923 \approx \frac{9}{13}
   \leq \tau^{(2)}_5 \leq 0.693943.
\]
\end{theorem}

\begin{proof}
The lower bounds were verified in Corollary \ref{lowerbounds}. We 
take the following general approach for finding an upper bound 
for $\tau^{(2)}_k$.  First, choose a specific value of $n$ (not 
too large), and then (carefully) choose a weighting function 
$\omega$ for $K[n,n]$.  Next, use binary linear programming 
(involving $n^2$ binary variables) to find 
$M_k(K[n,n],\omega)$.  This yields the upper bound
\[ \tau_k^{(2)} \leq \frac{M_k(K[n,n],\omega)}{W(\omega)}, \]
as given by Lemma \ref{weightedkings}.

For the weighting function $\omega_1$ shown (on the left) in 
Table \ref{weight4-5} for $K[10,10]$ with total weight 
$W(\omega_1)=280$, we computed $M_4(K[10,10],\omega_1)=171$.  
Similarly, the weighting function $\omega_2$ (on the right) in 
Table \ref{weight4-5} for $K[11,11]$ with $W(\omega_2)=2656$ gave 
us $M_5(K[11,11],\omega_2)=1844$.  These calculations prove 
the upper bounds
\[\tau^{(2)}_4 \leq \frac{171}{280} \approx 0.61071
  \quad\mbox{and}\quad
  \tau^{(2)}_5 \leq \frac{1844}{2656} = \frac{461}{664} \approx 
  0.69428.
\]
The tighter bounds stated in the theorem for $\tau^{(2)}_4$ and
$\tau^{(2)}_5$ were obtained by using  
significantly more complicated weighting functions on $K[12,12]$ 
and $K[13,13]$, respectively. These 
weighting functions have been posted on the Web \cite{website}.
\end{proof}


\begin{table}
\scriptsize$
  \begin{tabular}{|c|c|c|c|c|c|c|c|c|c|}
    \hline \strut 0& 1& 1& 1& 1& 1& 1& 1& 1& 0\\
    \hline 1& 2& 3& 3& 2& 2& 3& 3& 2& 1\\
    \hline 1& 3& 5& 4& 4& 4& 4& 5& 3& 1\\
    \hline 1& 3& 4& 6& 5& 5& 6& 4& 3& 1\\
    \hline 1& 2& 4& 5& 7& 7& 5& 4& 2& 1\\
    \hline 1& 2& 4& 5& 7& 7& 5& 4& 2& 1\\
    \hline 1& 3& 4& 6& 5& 5& 6& 4& 3& 1\\
    \hline 1& 3& 5& 4& 4& 4& 4& 5& 3& 1\\
    \hline 1& 2& 3& 3& 2& 2& 3& 3& 2& 1\\
    \hline 0& 1& 1& 1& 1& 1& 1& 1& 1& 0\\
    \hline
  \end{tabular}
\quad
  \begin{tabular}{|c|c|c|c|c|c|c|c|c|c|c|}
    \hline 2& 4& 6& 6& 6& 6& 6& 6& 6& 4& 2\\
    \hline 4& 11& 15& 17& 20& 22& 20& 17& 15& 11& 4\\
    \hline 6& 15& 23& 27& 36& 40& 36& 27& 23& 15& 6\\
    \hline 6& 17& 27& 36& 45& 59& 45& 36& 27& 17& 6\\
    \hline 6& 20& 36& 45& 49& 49& 49& 45& 36& 20& 6\\
    \hline 6& 22& 40& 50& 49& 48& 49& 50& 40& 22& 6\\
    \hline 6& 20& 36& 45& 49& 49& 49& 45& 36& 20& 6\\
    \hline 6& 17& 27& 36& 45& 59& 45& 36& 27& 17& 6\\
    \hline 6& 15& 23& 27& 36& 40& 36& 27& 23& 15& 6\\
    \hline 4& 11& 15& 17& 20& 22& 20& 17& 15& 11& 4\\
    \hline 2& 4& 6& 6& 6& 6& 6& 6& 6& 4& 2\\
    \hline
  \end{tabular}
$
\caption{
Weighting functions $\omega_1$ for $k=4$ and $\omega_2$ or $k=5$.
}
\label{weight4-5}
\end{table}

In the above proof, we say that the weighting function should be 
chosen carefully.  In fact, the weighting functions $\omega_1$ 
and $\omega_2$ were chosen to yield the best possible upper 
bounds for the particular kings graphs considered ($K[10,10]$ and 
$K[11,11]$, respectively).  We close this section with a quick 
explanation of how these were found.

Suppose we have a finite collection $\mathcal{S}$ of 
$k$-dependent subsets of $V(K[n,n])$. Consider the problem of  
minimizing a scalar $\theta$ over all pairs $(\theta,\omega)$ 
subject to the constraints
\begin{gather*}
  \omega(\vx)\geq0,\; \forall \vx\in V(K[n,n]),\\
  \sum_{\vx\in V(G)}\omega(\vx)=1,\\
  \sum_{\vx\in V(G)}\omega(\vx) \cdot x_\vx\leq\theta,
    \;\forall S\in\mathcal{S},
\end{gather*}
where we identify a set $S\in\mathcal{S}$ with its characteristic 
vector $(x_v:v\in V(K[n,n]))$.  The above optimization problem is 
a continuous-variable LP that can be solved in just seconds even 
when $|\mathcal{S}|\approx100,000$, provided that $n<50$.  Note 
that if the collection $\mathcal{S}$ contained all $k$-dependent 
subsets of $V(K[n,n])$, then the solution of this continuous LP 
would satisfy $\theta=M_k(K[n,n],\omega)$.  Moreover, this 
$\theta$ would be the best possible upper bound using $K[n,n]$ 
under Lemma \ref{weightedkings}.  However, $\mathcal{S}$ can also 
have this property and be considerably smaller than the entire 
collection of $k$-dependent sets. To check if $\mathcal{S}$ is 
sufficient for this purpose, solve the binary LP to maximize the 
$\omega$-weight, then obtain the value $M_k(K[n,n],\omega)$ and a 
corresponding maximum $\omega$-weight $k$-dependent set $S$.  If 
the value of $M_k(K[n,n],\omega)$ equals the optimal $\theta$ 
from the continuous LP, then the collection $\mathcal{S}$ is 
sufficient and we are done; otherwise, replace $\mathcal{S}$ by 
$\mathcal{S}\cup\{S\}$ and solve the above continuous LP again.  
This procedure is necessarily finite and guaranteed to find the 
best bound; even if terminated early it can find a very good 
bound, such as that provided by $\omega_3$ in the proof of 
Theorem \ref{maintheorem}.

The linear programming techniques for $k$-dependence and related 
problems lend themselves nicely toward undergraduate and Master's 
level research projects, provided the research supervisor can 
help with the details of getting LP packages to perform well.  In 
particular, there are many opportunities for using Lemma 
\ref{weightedtransitive} to do further research.  Any Cayley 
graph, for instance, is vertex-transitive, so $k$-dependence 
numbers for Cayley graphs are amenable to study in this manner.  


\section{Other Techniques: $\tau_1^{(d)}$ and $\tau_3^{(2)}$}
\label{others}

The binary LP technique used in Theorem \ref{maintheorem} also 
improves the upper bounds on $\tau_1^{(2)}$ and $\tau_3^{(2)}$, 
but is not nearly as effective in these two cases.  Fortunately, 
exact values for each can be found by other means.

The value $\tau_1^{(2)}=1/3$ is a special case of the next result.

\begin{theorem}\label{k=1} For any dimension $d$, we
have $\tau^{(d)}_1=2^{2-d}3^{-1}$.
\end{theorem} 

\begin{proof}
For any vertex $\vx\in V(T^{(d)}[n])$, define
\[ B_\vx=
   \left\{\yx\in N(\vx)\;\left|\; y_i-v_i\equiv 0\mbox{ or }1
       \hskip-.6em\pmod{n}
   \right.\right\}.
\]
Clearly, $|B_\vx|=2^d$. Also, $B_\vx\cap B_\ux\not=\emptyset$ 
precisely when $\vx$ and $\ux$ are neighbors. Now consider any 
$1$-dependent set $S\subseteq V(T^{(d)}[n])$ of maximum 
cardinality.  Suppose that $\vx,\ux\in S$ are neighbors for which
$v_i-u_i\equiv1\pmod{n}$ for some $i$. If $\yx\in B_\vx\cap 
B_\ux$, then $w_i=v_i$ and $|B_\vx\cap B_\ux|\leq2^{d-1}$. 
Because $S$ is 1-dependent, there are at most $|S|/2$ pairs 
$(\vx,\ux)$ for which $B_\vx$ and $B_\ux$ intersect.  Therefore
\[ |T^{(d)}[n]|
   \geq \left|\bigcup_{\vx\in S}B_\vx\right|
   \geq |S|\cdot2^d - \frac{|S|}{2}2^{d-1}
   =|S|\cdot 2^{d-2}3,
\]
so $\tau_1^{(d)}\leq 2^{d-2}3^{-1}$.
For the reverse inequality, consider $G=T[n_1,\ldots,n_d]$ with 
$n_1$ divisible by 3 and $n_i$ even for each $i>1$. Then the 
$1$-dependent set
\[ 
   \left\{\vx \in V(G) \;\left|\;
      \mbox{$v_1 \not\equiv 0\hskip-.6em\pmod{3}$ and 
       $v_i \equiv 0\hskip-.6em\pmod{2},\, \forall i>1$}
   \right.\right\}
\]
has cardinality $|V(G)|/(2^{d-2}3)$, so 
$\tau_1^{(d)}\geq\tau_1(T[n_1,\ldots,n_d])\geq2^{2-d}3^{-1}$.
\end{proof}

We close by deriving the exact value $\tau^{(2)}_3=1/2$, which is 
surprising in that we cannot improve upon the 2-dependent density 
by admitting a third neighboring king.  Our proof uses a 
``taxation'' argument, like those in \cite{ffpp}.  

\begin{theorem}
In $T[m,n]$, every 3-dependent set has at least $mn/2$ vertices.  
Therefore, $\tau_3^{(2)} = 1/2$.
\end{theorem}

\begin{proof}[Sketch of proof]
Consider a $3$-dependent set $S$ of vertices in $T[m,n]$. The 
idea of taxation is to start with \$1 at every vertex \emph{not} 
in $S$, and then redistribute those funds in such a way that each 
member of $S$ receives at least \$1.  After redistribution, the 
vertices in $S$ collectively share a total of at least $|S|$ 
dollars, whereas the same total cannot exceed the 
$|V(T[m,n])\setminus S|$ dollars originally distributed over the 
vertices in the complement of $S$.  Therefore 
$|V(T[m,n])|-|S|\geq |S|$, and so $mn/2\geq|S|$. The theorem then 
follows by observing that $\tau_3^{(2)} \geq \tau_2^{(2)} =1/2$. 

A taxation argument hinges on the particular rule used to 
redistribute funds.  To describe a suitable rule, we view the 
vertices as squares on a toroidal chessboard.  For each vertex 
$\vx$ of $T[m,n]$, let $N_{\text{side}}(\vx)$ denote the 
neighbors of $\vx$ that share a side with $\vx$ and let 
$N_{\text{corner}}(\vx)$ denote the neighbors that do not share a 
side with $\vx$.
Let $r_0(\vx)$ denote the amount of money initially available at 
vertex $\vx$, so that $r_0(\vx)=1$ for $\vx\not\in S$ and 
$r_0(\vx)=0$ for $\vx\in S$; let $t(\vx)=1-r_0(\vx)$ denote the 
``target'' amount for $\vx$. We redistribute the money by the 
following three steps, where $r_i(\vx)$ denotes the amount at 
vertex $\vx$ immediately after step $i$:   
\begin{itemize}
\item[1.] Each $\vx$ distributes its ``surplus'' $\max\{r_0(\vx)-t(\vx),0\}$
evenly among those $\ux\in N_{\text{side}}(\vx)$ for which 
$r_0(\ux)<t(\ux)$, if any.
\item[2.] Each $\vx$ with $r_1(\vx)>t(\vx)$ transfers
the amount $\max\{t(\ux)-r_1(\ux),0\}$ to each $\ux \in 
N_{\text{corner}}(\vx)$.
\item[3.] Each $\vx$ transfers the amount $\max\{r_2(\vx)-t(\vx),0\}$ to each 
$\ux\in N_{\text{side}}(\vx)$ for which $r_2(\ux)<t(\ux)$.    
\end{itemize}
\noindent A case-by-case examination of the neighborhood 
possibilities for members of the 3-dependent set $S$ verifies 
that $r_i(\vx)\geq0$ and $r_3(\vx)\geq t(\vx)$ for all $\vx$.
\end{proof}


\bigskip

\noindent \textbf{\large Acknowledgements.}\ \ Drs.~David 
Woolbright and Tim Howard at Columbus State University first 
brought our attention to this problem.  Some examples on large 
chessboards ($n>20$) were constructed using a program written by 
Mike McCoy.  Linear programming computations were carried out 
using the LPsolve and ILOG Cplex software packages, with the most 
difficult calculations performed on Miami University's 
``RedHawk'' computing cluster.

\end{document}